\documentclass[12pt]{article}

\newcommand{\be}{\begin{equation}}
\newcommand{\ee}{\end{equation}}

\begin{document}

\begin{center}
Journal of Physics A. Vol.38. No.10/11. (2005) pp.2145-2155. 
\end{center}

{\Large \bf Phase-Space Metric for Non-Hamiltonian Systems}

\vskip 7mm
\begin{center}
{\large \bf Vasily E. Tarasov } \\
{\it Skobeltsyn Institute of Nuclear Physics, 
Moscow State University, Moscow 119992, Russia \\
E-mail: tarasov@theory.sinp.msu.ru } \\
\end{center}

\begin{abstract}
We consider an invariant skew-symmetric phase-space metric 
for non-Hamiltonian systems. 
We say that the metric is an invariant if the metric tensor field 
is an integral of motion. 
We derive the time-dependent skew-symmetric phase-space metric 
that satisfies the Jacobi identity. 
The example of non-Hamiltonian systems with linear friction term
is considered. 
\end{abstract}


\section{Introduction}

The dynamics of Hamiltonian systems is characterized by conservation 
of  phase-space volume under time evolution. 
This conservation of the phase volume is a cornerstone of conventional 
statistical mechanics of Hamiltonian systems. 
At a mathematical level, conservation of phase-space volume is 
considered as a consequence of the existence of an 
invariant symplectic form (skew-symmetric phase-space metric) 
in the phase-space of Hamiltonian systems \cite{Godb,DNF,Fom}. 

The classical statistical mechanics 
of non-Hamiltonian systems is of strong theoretical interest
\cite{Tuck1,Tuck2,Tuck3,Ram1,Ram2,Sergi1,Sergi2,Sergi3,Tarmpl,chaos,Ezra1,Ezra2,Tarpre,Tarpla}.
Non-Hamiltonian systems have been used in molecular dynamics 
simulation to achieve the calculation of statistical averages 
in various ensemble \cite{Tuck2,Tuck3,Tarmpl,Nose}, 
and in the treatment of 
nonequilibrium  steady states \cite{Tarpre,Tarpla,Evans,Hoover}. 
Non-Hamiltonian systems are characterized by nonzero phase 
space compressibility, and 
the usual phase-space volume is no longer necessarily conserved.

Tuckerman {\it et al.} have argued \cite{Tuck2,Tuck3} that there 
is a measure conservation law that involves 
a nontrivial phase-space metric. This suggests that phase-space 
should be carefully treated using 
the general rules of the geometry of manifolds \cite{Godb,DNF}. 
Tuckerman et al. have applied the concepts of  
Riemannian geometry to the classical statistical mechanics 
of non-Hamiltonian systems \cite{Tuck1,Tuck2,Tuck3}. 
Tuckerman {\it et al.} have argued that, through 
introduction of metric determinant factors $\sqrt{g({\bf x},t)}$,
it is possible to define an invariant phase-space measure 
for non-Hamiltonian systems.
In their approach the metric determinant factor $\sqrt{g({\bf x},t)}$,  
where $g({\bf x},t)$ is the determinant of the metric tensor, 
is defined by the compressibility of non-Hamiltonian systems.
However Tuckerman {\it et al.} consider only the 
determinant $g({\bf x},t)$ of the metric. 
The phase-space metric is not considered in \cite{Tuck1,Tuck2,Tuck3}.  
Note that Tuckerman {\it et al.} suppose that the metric determinant factor
is connected with symmetric phase-space metric.
It can be proved that the proposal to use an invariant 
time-dependent metric determinant factor
in the volume element corresponds precisely 
to finding a skew-symmetric phase-space metric (symplectic form)
that is an integral of motion. 
Therefore we must consider the skew-symmetric phase-space metric.

Sergi \cite{Sergi2,Sergi3} has considered an antisymmetric 
phase-space tensor field, whose elements are 
general function of phase-space coordinates. 
In \cite{Sergi2,Sergi3}, the generalization of Poisson brackets 
for the non-Hamiltonian systems was suggested. 
However the Jacobi identity is not satisfied by the generalized brackets
and skew-symmetric phase-space metric. 
As a result the algebra of phase-space functions is not 
time translation invariant. The generalized brackets do not define 
a Lie algebra in phase-space. Note that the generalized brackets of 
two constant of motion is no longer a constant of motion.

In the present paper we consider an invariant 
skew-symmetric (antisymmetric) phase-space metric for non-Hamiltonian systems. 
We say that the metric is an invariant, if the metric tensor field 
is an integral of motion. 
We define the phase-space metric such that 
the Jacobi identity is satisfied.
The suggested skew-symmetric phase-space metric allows us
to introduce the generalization of the Poisson brackets for 
non-Hamiltonian systems such that 
the Jacobi identity is satisfied by the generalized Poisson brackets.
As a result the algebra of phase-space functions is time 
translation invariant. The generalized Poisson brackets define 
a Lie algebra in phase-space. The suggested Poisson brackets of 
two constant of motion is a constant of motion.

In section 2, the definitions of the antisymmetric phase-space metric, 
mathematical background and notations are considered.
In section 3, we define the non-Hamiltonian systems, and 
consider the Helmholtz conditions. 
In section 4, we consider the time evolution of phase-space metric.
We derive the phase-space metric that is an integral of motion.
In section 5, the generalized Poisson brackets for non-Hamiltonian
systems are defined. 
In section 6, the example of phase-space metric for non-Hamiltonian system 
with the linear friction term is considered. 
Finally, a short conclusion is given in section 7.

\section{Phase-space metric}

The 2n-dimensional differentiable manifold is denoted $M$. 
Coordinates are ${\bf x}=(x^1,...,x^{2n})$. 
We assume the existence of a time-dependent metric tensor field
$\omega_{kl}({\bf x},t)$ on the manifold $M$.  
We can define a differential 2-form
\be \label{om} \omega=\omega_{kl}({\bf x},t) dx^k \wedge dx^l , \ee
where $\omega_{kl}=\omega_{kl}({\bf x},t)$ is a skew-symmetric tensor
$\omega_{kl}=-\omega_{lk}$, 
and the tensor elements $\omega_{kl}({\bf x},t)$ are explicit
functions of time. 
Here and later we mean the sum on the repeated index
$k$ and $l$ from 1 to 2n.

We suppose that the differential 2-form $\omega$ is a
closed nondegenerated form: \\

\noindent
1) If the metric determinant is not equal to zero
\be \label{g} g({\bf x},t)=det(\omega_{kl}({\bf x},t))\not =0 \ee
for all points ${\bf x}\in M$, then the form $\omega$ is nondegenerated. \\
2) If the Jacobi identity 
\be \label{Ji} \partial_k \omega_{lm}+\partial_l \omega_{mk}+\partial_m \omega_{kl}=0 , 
\quad \partial_k= \partial / \partial x^k \ee
for the metric $\omega_{kl}=\omega_{kl}({\bf x},t)$ is satisfied, 
then the differential 2-form $\omega$ is closed ($d\omega=0$).
 
Phase-space is therefore assumed to be a symplectic manifold. \\

{\bf Definition 1}. 
{\it A {\bf symplectic manifold} is a differentiable manifold $M$
with a closed nondegenerated differential 2-form $\omega$.} \\

For symplectic manifold, we have the phase-space volume element
\[ v=\frac{1}{n!} \omega^n=\frac{1}{n!}\omega \wedge ... \wedge \omega , \]
This differentiable $2n$-form can be represented by 
\[ v=\sqrt{g({\bf x},t)} dx^1 \wedge ... \wedge dx^{2n} , \]
where $g({\bf x},t)$ is defined by equation (\ref{g}).
The nondegenerated condition ($g({\bf x},t)\not=0$)
for the metric $\omega_{kl}({\bf x},t)$ is equivalent to 
the condition $\omega^n\not=0$ or $v\not=0$.

It is known \cite{DNF} that there exists the local coordinates $(q,p)$ 
such that
\be \label{oqp} \omega=\delta_{ij} dq^i \wedge dp^j . \ee
Here and later we mean the sum on the repeated index
$i$ and $j$ from 1 to $n$.

\section{Non-Hamiltonian system}

The dynamics is described by a smooth vector field 
${\bf X}={\bf X}({\bf x})$,
\be \label{eq1} \frac{d{\bf x}}{dt}={\bf X} \ee
with components $X^k$ in basis $\partial_k=\partial / \partial x^k$.
For simplicity, we consider the case where the vector 
field ${\bf X}$ is time independent.
In local coordinates $\{x^k\}$, equation (\ref{eq1}) has the form
\be \label{eq2} \frac{dx^k}{dt}=X^k. \ee

Consider now the definition of the Hamiltonian systems, 
which is used in \cite{Godb}. \\

{\bf Definition 2}. 
{\it A classical system (\ref{eq1}) on the symplectic manifold 
$(M,\omega)$ is called a {\bf Hamiltonian system} if the differential 
1-form $\omega(X)$ is a closed form
\[ d \omega(X)=0 , \]
where $\omega(X)=i_X \omega$ is the contraction (interior product)
of the 2-form $\omega$ with vector ${\bf X}$, and
$d$ is the exterior derivative. 

A classical system (\ref{eq1}) on the symplectic manifold $(M,\omega)$
is called a {\bf non-Hamiltonian system} if the 
differential 1-form  $\omega(X)$ is nonclosed $d \omega(X)\not=0$}. \\

{\bf Proposition 1}.
{\it The classical system (\ref{eq1}) is a Hamiltonian system if 
the conditions
\be \label{eq0} J_{kl}(\omega,{\bf x},t)\equiv
\partial_k (\omega_{lm} X^m)-
\partial_l (\omega_{km} X^m) =0 \ee
are satisfied.} \\

{\bf Proof}.  
In the local coordinates $\{x^k\}$, we have
\[ \omega(X)=X_kdx^k=\omega_{kl} X^l dx^k,  \]
where $X_k=\omega_{kl}X^l$. 
In this case, the exterior derivative of 1-form $\omega(X)$ is
\[ d\omega(X)=d(X_kdx^k)=\partial_l X_k dx^l \wedge dx^k . \] 
Using $da \wedge db=- db \wedge da$, we get
\[ d \omega(X)=\frac{1}{2}(\partial_k X_l-\partial_l X_k) dx^k \wedge dx^l . \]
As the result we have the differential 2-form
\be \label{Jkl} d \omega(X)=
\frac{1}{2} \Bigl( \partial_k (\omega_{lm} X^m)-
\partial_l (\omega_{km} X^m) \Bigr) 
 dx^k \wedge dx^l.  \ee
This differential 2-form is a symplectic form,
which can be called "non-Hamiltonian symplectic form". 
If the Helmholtz conditions (\ref{eq0})
are satisfied, then the differential 1-form $\omega(X)$ 
is closed ($d\omega(X)=0$), and 
the classical system (\ref{eq1}) is a Hamiltonian system.  \\

Let us consider the canonical coordinates 
${\bf x}=(x^1,...,x^n,x^{n+1},...,x^{2n})=(q^1,...,q^n,p^{1},...,p^{n})$. 
Equation (\ref{eq2}) can be written as
\be \label{em1} \frac{dq^i}{dt}=G^i(q,p), 
\quad \frac{dp^i}{dt}=F^i(q,p) . \ee

{\bf Corollary}. 
{\it If the right-hand sides of equations (\ref{em1})
satisfy the Helmholtz conditions \cite{Helm,Tartmf3} for the 
phase-space with (\ref{oqp}), which have the following form 
\be \label{HC1} \frac{\partial G^{i}}{\partial p^j}-
\frac{\partial G^{j}}{\partial p^i}= 0, \ee
\be \label{HC2} \frac{\partial G^{j}}{\partial q^i}+
\frac{\partial F^{i}}{\partial p^j}=0, \ee
\be \label{HC3} \frac{\partial F^{i}}{\partial q^j}-
\frac{\partial F^{j}}{\partial q^i}= 0, \ee
then classical system (\ref{em1}) is a Hamiltonian system.} \\

{\bf Proof}. 
In the canonical coordinates $(q,p)$, the vector field ${\bf X}$
has the components $(G^i,F^i)$, which are used in equation (\ref{em1}). 
The 1-form $\omega(X)$ is defined by the following equation:
\[ \omega(X)=\frac{1}{2}(G_idp^i-F_idq^i), \]
where $G_i=\delta_{ij}G^j$ and $F_i=\delta_{ij}F^j$. 
The exterior derivative for this form can now be written by the relation
\[ d\omega(X)=\frac{1}{2}\Bigl( d(G_idp^i)-d(F_idq^i)\Bigr).  \]
It now follows that
\[ d\omega(X)= \frac{1}{2}\left( \frac{\partial G_i}{\partial q^j} 
dq^j \wedge dp^i
+\frac{\partial G_i}{\partial p^j} dp^j \wedge dp^i-
\frac{\partial F_i}{\partial q^j} dq^j \wedge dq^i-
\frac{\partial F_i}{\partial p^j} dp^j \wedge dq^i \right) . \]
This equation can be rewritten in an equivalent form
\[ d\omega(X)=
\frac{1}{2}\left( \frac{\partial G_j}{\partial q^i} 
+\frac{\partial F_i}{\partial p^j}\right) dq^i \wedge dp^j+
\frac{1}{4}\left( \frac{\partial G_j}{\partial p^i}
-\frac{\partial G_i}{\partial p^j} \right)dp^i \wedge dp^j
+\frac{1}{4}\left(\frac{\partial F_i}{\partial q^j}-
\frac{\partial F_j}{\partial q^i} \right) dq^i \wedge dq^j .  \]
Here we use the skew-symmetry of $dq^i \wedge dq^j$ and $dp^i \wedge dp^j$
with respect index $i$ and $j$.
It is obvious that conditions (\ref{HC1}) - (\ref{HC3}) 
lead to the equation $d\omega(X)=0$.

\section{Time evolution of phase-space metric}

Let us find a time-dependent symplectic 2-form $\omega$ that 
satisfies the equation $d \omega / dt=0$. 

It is known \cite{DNF,Fom} as the following proposition. \\

{\bf Proposition 2}.
{\it If the system $\dot{\bf x}={\bf X}$ 
on the symplectic manifold $(M,\omega)$ 
with time-independent symplectic form 
($\partial \omega_{kl}/ \partial t=0$)
is a Hamiltonian system, then 
differential 2-form $\omega$ is conserved, i.e.,
$d \omega /dt=0$.} \\

{\bf Proof}. 
The proof of this theorem is considered in \cite{DNF,Fom}. \\

Let us consider a generalization of this proposition. \\

{\bf Proposition 3}. 
{\it If the time-dependent metric $\omega_{kl}=\omega_{kl}({\bf x},t)$ 
is a skew-symmetric 
metric ($\omega_{kl}=-\omega_{lk}$) that is satisfied by 
the Jacobi identity (\ref{Ji}), 
and the system is defined by equation (\ref{eq2}), 
then the total time derivative of the differential 2-form (\ref{om})
is given by }
\be \label{do} \frac{d \omega}{dt} =
\Bigl(\frac{\partial \omega_{kl}}{\partial t}
-  \partial_k (\omega_{lm} X^m)+
\partial_l (\omega_{km} X^m) \Bigr)dx^k \wedge dx^l . \ee

{\bf Proof}. 
The time-derivative of the time-dependent symplectic form $\omega$ is given by
\[ \frac{d\omega}{dt}=
\frac{d}{dt}\Bigl(\omega_{kl}({\bf x},t) dx^k \wedge dx^l\Bigr)=
\frac{d\omega_{kl}}{dt} dx^k \wedge dx^l
+\omega_{kl} d \left(\frac{dx^k}{dt}\right) \wedge dx^l
+\omega_{kl} dx^k \wedge d\left(\frac{dx^l}{dt}\right). \]
Then, using the equation
\[ \frac{d \omega_{kl}}{dt}=\frac{\partial \omega_{kl}}{\partial t}+
\frac{\partial \omega_{kl}}{\partial x^m} \frac{dx^m}{dt} =
\frac{\partial \omega_{kl}}{\partial t}+
X^m \partial_m \omega_{kl}, \]
and equation (\ref{eq2}), we find that
\[ \frac{d\omega}{dt}=
\left( \frac{\partial \omega_{kl}}{\partial t}+
X^m \partial_m \omega_{kl}\right) dx^k \wedge dx^l
+\omega_{kl} d X^k \wedge dx^l
+\omega_{kl} dx^k \wedge d X^l . \]
Using $dX^k=\partial_m X^k dx^m$, we have
\[ \frac{d\omega}{dt}=
\left( \frac{\partial \omega_{kl}}{\partial t}+
X^m \partial_m \omega_{kl}\right) dx^k \wedge dx^l
+\omega_{kl} \partial_m X^k dx^m \wedge dx^l
+\omega_{kl}  \partial_m X^l dx^k \wedge d x^m . \]
This expression can be rewritten in an equivalent form
\[ \frac{d\omega}{dt}=\left( \frac{\partial \omega_{kl}}{\partial t}+
X^m \partial_m \omega_{kl}+
\omega_{ml} \partial_k X^m +\omega_{km} \partial_l X^m 
\right) dx^k \wedge d x^l . \]
Using the rule of term-by-term differentiation in the form
\[ \omega_{ml} \partial_k X^m=
\partial_k(\omega_{ml} X^m)-X^m \partial_k \omega_{ml}; 
\quad  \omega_{km} \partial_l X^m=
\partial_l(\omega_{km} X^m)-X^m \partial_l \omega_{km},  \]
we get the following equation:
\[ \frac{d\omega}{dt}=\left( \frac{\partial \omega_{kl}}{\partial t}+
X^m (\partial_m \omega_{kl}- \partial_k\omega_{ml}-\partial_l \omega_{km})
+\partial_k(\omega_{ml} X^m)+\partial_l(\omega_{km} X^m)
\right) dx^k \wedge d x^l . \]
Using the Jacobi identity (\ref{Ji}), and skew symmetry 
$\omega_{ml}=-\omega_{lm}$, $\omega_{km}=-\omega_{mk}$, we have 
\[ \frac{d\omega}{dt}=\left( \frac{\partial \omega_{kl}}{\partial t}
-\partial_k(\omega_{lm} X^m)+\partial_l(\omega_{km} X^m)
\right) dx^k \wedge d x^l . \]
As the result, we obtain equation (\ref{do}) for 
the total time derivative of symplectic form. \\

The total time derivative of the symplectic form 
is defined by equation (\ref{do}). 
If the total derivative is zero, than 
we have the integral of motion or invariant. 
It is easy to see that
the  differentiable 2-form $\omega$ is invariant if
the phase-space metric $\omega_{kl}({\bf x},t)$ 
is satisfied by the equation
\be \label{7} \frac{\partial \omega_{kl}}{\partial t}
= \partial_k (\omega_{lm} X^m)-\partial_l (\omega_{km} X^m) . \ee
This equation can be rewritten in an equivalent form
\[ \frac{\partial \omega_{kl}}{\partial t}=\hat J^{ms}_{kl} \omega_{ms} ,\]
where the operator $\hat J$ is defined by the equation 
\be \label{hJ} \hat J^{ms}_{kl}=\frac{1}{2}\Bigl(
(\delta^m_l\partial_k-\delta^m_k \partial_l)X^s -
(\delta^s_l\partial_k-\delta^s_k \partial_l)X^m \Bigr) . \ee

\vskip 3mm

{\bf Proposition 4}. 
{\it The  differentiable 2-form $\omega$ is invariant 
(is an integral of motion for non-Hamiltonian system (\ref{eq2}))
if the phase-space metric $\omega_{kl}({\bf x},t)$ 
is defined by the equation
\be \label{o-sol} \omega_{kl}({\bf x},t)=\Bigl(exp (t \ \hat J)\Bigr)^{ms}_{kl}
\omega_{ms}({\bf x},0) . \ee 
Here $\hat J$ is an operator that is defined by equation (\ref{hJ}). } \\

{\bf Proof}. 
Let us consider the formal solution of equation (\ref{7}) in the form
\[ \omega_{kl}({\bf x},t)=\sum^{\infty}_{n=1}
\frac{t^n}{n!} \omega^{(n)}_{kl}({\bf x}) , \]
where $\omega^{(n)}_{kl}=-\omega^{(n)}_{lk}$.  
In this case, the time-independent tensor fields 
$\omega^{(n)}_{kl}({\bf x})$ are defined by the recursion relation
\[ \omega^{(n+1)}_{kl}=\partial_k (\omega^{(n)}_{lm} X^m)-
\partial_l (\omega^{(n)}_{km} X^m) . \]
This equation can be rewritten in an equivalent form
\[ \omega^{(n+1)}_{kl}=
(\delta^m_l\partial_k-\delta^m_k \partial_l)(X^s\omega^{(n)}_{ms}) . \]
Using the skew symmetry of the $\omega^{(n)}_{ns}$,  we have
\[ \omega^{(n+1)}_{kl}=
\frac{1}{2} \Bigl(
(\delta^m_l\partial_k-\delta^m_k \partial_l)X^s -
(\delta^s_l\partial_k-\delta^s_k \partial_l)X^m 
\Bigr) \omega^{(n)}_{ms} . \]
This relation can be represented in the form
\[ \omega^{(n+1)}_{kl}=\hat J^{ms}_{kl} \omega^{(n)}_{ms} ,\]
where the operator $\hat J$ is defined by equation (\ref{hJ}). 
Therefore the invariant phase-space metric is 
defined by the following equation: 
\[ \omega_{kl}({\bf x},t)=\sum^{\infty}_{n=0} \frac{t^n}{n!} 
\Bigl({\hat J}^n \Bigr)^{ms}_{kl} \omega_{ms}({\bf x},0)=
\Bigl(exp (t \ \hat J)\Bigr)^{ms}_{kl}
\omega_{ms}({\bf x},0) . \]
As the result  we have equation (\ref{o-sol}).  \\

{\bf Comments}. 
We prove that equation (\ref{7}) for the elements of phase-space metric
can be expressed in terms of an operator $\hat J$ such that
\be \label{7-2} \frac{\partial \omega_{kl}}{\partial t}=
\hat J^{ms}_{kl} \omega_{ms}, \ee
where $\hat J^{ab}_{kl}$ is defined by equation (\ref{hJ}). 
In this case, the time evolution of the phase-space metric
from initial condition $\omega_{kl}({\bf x}_0,0)$ to a value
$\omega_{kl}({\bf x},t)$ at the time $t$ can be written 
by the equation
\be \label{o-sol2} \omega_{kl}({\bf x},t)=
\Bigl(exp (t \ \hat J)\Bigr)^{ms}_{kl}
\omega_{ms}({\bf x},0) . \ee 
Here the matrix exponential operator, $exp(t \hat J )$,  
can be called the metric propagator. 
The operator $\hat J$ can be considered as a metric analog of the Liouville
operator. 
Note that we can use the canonical coordinates $(q,p)$ for $t=0$ and
the coordinate-independent initial metric: 
$\omega_{kl}({\bf x},0)=\omega^{(0)}_{kl}=const$, 
$\partial_m \omega_{kl}({\bf x},0)=0$. 
Introducing time step $\Delta t=t/N$, we get
\[ \omega_{kl}({\bf x},t)=\Bigl(exp (t \ \hat J)\Bigr)^{ms}_{kl}
\omega_{ms}({\bf x},0)=\Bigl([exp (\Delta t \ \hat J)]^N\Bigr)^{ms}_{kl}
\omega_{ms}({\bf x},0) . \]
It is natural to approximate the short-time propagator $exp(t\hat J)$
using the Trotter theorem \cite{T1,T2,T3}.
The vector field ${\bf X}={\bf X(x)}$ can be represented in the form
\[ {\bf X}={\bf X}_1+{\bf X}_2 , \]
where ${\bf X}_1$ is a Hamiltonian term such that
\[ X^k_1=\omega^{kl} \frac{\partial H}{\partial x^l} , \] 
and ${\bf X}_2$ is a friction (non-Hamiltonian) term. 
As the result we have
\be \label{hJ12} (\hat J_{1;2})^{ms}_{kl}=\frac{1}{2}\Bigl(
(\delta^m_l\partial_k-\delta^m_k \partial_l)X^s_{1;2} -
(\delta^s_l\partial_k-\delta^s_k \partial_l)X^m_{1;2} \Bigr) . \ee
An important consequence would be the ability to formulate rigorous 
numerical integration algorithms based on Trotter-type splittings 
\cite{T1,T2,T3} of the classical propagator $exp(\Delta t \hat J)$. 
The metric propagator for the small time step is
\[ exp[\Delta t (\hat J_1+\hat J_2)]=
exp[\frac{1}{2}\Delta t \hat J_2] exp[\Delta t \hat J_1]
exp[\frac{1}{2}\Delta t \hat J_2]+O(\Delta t^3) . \]
Finally, we obtain
\[ exp[t \hat J ]=exp[t (\hat J_1+\hat J_2)]=
\Bigl( exp[\frac{1}{2}\Delta t \hat J_2] exp[\Delta t \hat J_1]
exp[\frac{1}{2}\Delta t \hat J_2]\Bigr)^N+O(\Delta t^3) . \]

\section{Poisson brackets for non-Hamiltonian systems}

Let us consider the skew-symmetric tensor field 
$\omega^{kl}=\omega^{kl}({\bf x},t)$ 
that is defined by the equations
\[ \omega^{kl}({\bf x},t) \omega_{lm}({\bf x},t)=
\omega^{lk}({\bf x},t) \omega_{ml}({\bf x},t)=\delta^k_l . \]
As the result this tensor field satisfies the Jacoby identity 
\[ \omega^{kl}\partial_l \omega^{ms}+
\omega^{ml}\partial_l \omega^{sk}+
\omega^{sl}\partial_l \omega^{km}=0 . \]
It follows from the Jacoby identity for $\omega_{kl}$.

In proposition 3, we suggest the time-dependent phase-space metric 
$\omega_{kl}({\bf x},t)$, which satisfies the Jacoby identity. 
As the result we have Lie algebra that is defined by
the following brackets:
\be \label{PB} 
\{A,B\}=\omega^{kl}({\bf x},t) \partial_k A \partial_l B . \ee
It is easy to prove that these brackets are Poisson brackets. 

In the general case, 
the rule of term-by-term differentiation with respect to time that
has the form
\be \label{LR} \frac{d}{dt}\{A,B\}=\{\dot{A},B\}+\{A,\dot{B}\} , \ee
where $\dot A=dA/dt$, is not valid for non-Hamiltonian systems. 
In the general case, we have
\[ \frac{d}{dt}\{A,B\}=\{\dot{A},B\}+\{A,\dot{B}\}+J(A,B) \]
where 
\[ J(A,B)=\omega^{kl}_{(1)}({\bf x},t)  \partial_k A \  
\partial_l B, \quad 
\omega^{kl}_{(1)}({\bf x},t)=\omega^{km} \omega^{ls} 
(\partial_s X_m-\partial_m X_s) .\]
Note that  time evolution of the Poisson brackets (\ref{PB})
for non-Hamiltonian systems can be considered as t-deformation 
\cite{Ger} of the Lie algebra in phase-space.

If we use the invariant phase-space metric, then the rule (\ref{LR})
is valid. 
As the result the suggested Poisson brackets (\ref{PB}) of 
two constants of motion is a constant of motion and 
rule (\ref{LR}) is satisfied.

\section{Example: System with linear friction}

Let us consider the non-Hamiltonian system that is described 
by the equations:
\be \label{ex1}
\frac{dq^i}{dt}=\frac{\partial H}{\partial p^i} , \quad 
\frac{dp^i}{dt}=-\frac{\partial H}{\partial q^i}-K^i_j(t) p^j , \ee
where $H=T(p)+U(q)$ and $i;j=1,...,n$. 
Here $T(p)$ is a kinetic energy, $U(q)$ is a potential energy. 
The term $-K^i_j(t)p^j$ is a friction term. 
Usually, this system is described by the phase-space metric $\omega_{kl}$
that has the form
\be \label{GG}
|| \omega_{kl}||=\left(
\begin{array}{cc}
0&G \\
-G^{T}&0 \\
\end{array}
\right),
\ee
where the matrix $G=||g_{ij}||$ is equal to indentity 
matrix $E=||\delta_{ij}||$. Here $G^T$ is transpose matrix for 
the matrix $G$. 
The symplectic form is defined by equation (\ref{oqp}) in the form
\[ \omega=2g_{ij}dq^i\wedge dp^j =\delta_{ij} dq^i\wedge dp^j , \]
where $g_{ij}=(1/2)\delta_{ij}$.
The phase-space compressibility 
\be \label{psc} \kappa=\sum^n_{i=1} 
\left(\frac{\partial \dot{q^i}}{\partial q^i}
+\frac{\partial \dot{p^i}}{\partial p^i} \right) \ee
of the system (\ref{ex1}) is defined by the spur of the matrix $K=||K^i_j||$
in the form 
\[ \kappa=-Sp||K^i_j|| . \]

\vskip 3mm

{\bf Proposition 5}. 
{\it The invariant phase-space metric for the classical system (\ref{ex1})
has the form
\be \label{||o||} || \omega_{kl}(t)||=\left(
\begin{array}{cc}
0&G(t) \\
-G^{T}(t)&0 \\
\end{array}
\right),
\ee
where the matrix $G$ is defined by }
\be \label{G(t)} G(t)
= G(t_0)exp \int^t_{t_0} ||K^i_j(\tau)|| d \tau . \ee

{\bf Proof}. 
Suppose that phase-space metric depends on time $t$. Therefore 
the matrix $G$ and elements $g_{kl}$ are the functions of the variable $t$:
$G=G(t)$, $g_{ij}=g_{ij}(t)$. 
Let us consider the total time derivative of the symplectic form
\[ \omega=2g_{ij}(t)dq^i\wedge dp^j . \]
As the result, we have
\[ \frac{d \omega}{dt}=2\Bigl(  
\frac{dg_{ij}}{dt}-g_{im}K^{m}_j \Bigr) dq^i\wedge dp^j . \]
In order to have the invariant phase-space metric ($d \omega /dt=0$), 
we use the following equation: 
\[ \frac{dg_{ij}}{dt}= g_{im}K^{m}_j . \]
The matrix $G=||g_{ij}||$ is satisfied the matrix equation
\be \label{GK} 
\frac{dG(t)}{dt}=G(t)K(t),
\ee
where $K(t)=||K^i_j(t)||$ is a matrix of friction coefficients.
Suppose that $G(t_0)=E$, where $E$ is the identity matrix.
The solution of equation (\ref{GK}) has the form
\be 
\label{Gt} 
G(t)=G(t_0)exp \int^t_{t_0} K(\tau) d \tau =exp \int^t_{t_0} K(\tau) d \tau .\ee
If the matrix $K$ is diagonal matrix with elements 
$K^i_j(t)=K_j(t)\delta^i_j$, then we have the matrix elements
\[ g_{ij}(t) =\delta_{ij} exp\int^t_{t_0} K_j(\tau) d \tau .\]
As the result, the invariant phase-space metric for system
(\ref{ex1}) is defined by equations (\ref{||o||}) and (\ref{G(t)}). 

The metric determinant factor \cite{Tuck1,Tuck2,Tuck3}
\[ \sqrt{g({\bf x},t)}=\sqrt{det ||\omega_{kl}({\bf x},t)||} \]
for the phase-space metric (\ref{GG}) is defined by the relation
$\sqrt{g}=\sqrt{det(G G^T)}$.  
Using $det \ (AB)=(det \ A) \ (det \ B)$ and $det A^T=det A$, we have 
the metric determinant factor in the form $\sqrt{g}=|det G|$.
If the matrix $G$ has the form (\ref{Gt}), then 
\[ \sqrt{g(t)}=|det \ exp \int^t_{t_0} K(\tau) d\tau | . \]
Using $det \ exp \ A=exp \ Sp \ A$, we get that 
the invariant phase-space metric has the determinant 
that is connected with the phase-space compressibility
\cite{Tuck1,Tuck2,Tuck3} by the equation  
\be
\sqrt{g(t)}=
exp \int^t_{t_0} Sp K(\tau) d \tau =exp-\int^t_{t_0} \kappa(\tau) d \tau,
\ee
where $\kappa$ is the phase-space compressibility (\ref{psc}). 

For example, the system
\[ \frac{dq_1}{dt}=\frac{p_1}{m} , \quad 
\frac{dp_1}{dt}=-\frac{\partial U(q)}{\partial q_1}-K_1 p_1 , \]
\[ \frac{dq_2}{dt}=\frac{p_2}{m} , \quad 
\frac{dp_2}{dt}=-\frac{\partial U(q)}{\partial q_2}-K_2 p_2 , \]
has the invariant phase-space metric $\omega_{kl}(t)$ in the form
\be \label{ex2}
|| \omega_{kl}(t)||=\left(
\begin{array}{cccc}
0&0&e^{K_1t}&0 \\
0&0&0&e^{K_2t} \\
-e^{K_1t}&0&0&0 \\
0&-e^{K_2t}&0&0 \\
\end{array} \right). \ee
In this case, the metric determinant factor is equal to 
\[ \sqrt{g(t)}=\sqrt{det ||\omega_{kl}||} =e^{(K_1+K_2)t} . \]

\section{Conclusion}

Tuckerman {\it et al.} \cite{Tuck1,Tuck2,Tuck3} suggest a
formulation of non-Hamiltonian statistical mechanics which 
uses the invariant phase-space measure. 
The invariant measure is connected with phase-space metric.
Tuckerman {\it et al.} consider the properties of 
only  the determinant $g({\bf x},t)$ of the metric. 
The phase-space metric is not considered in \cite{Tuck1,Tuck2,Tuck3}.  
We consider the invariant phase-space metric for non-Hamiltonian systems.  
The proposal to use an invariant time-dependent metric determinant 
factor in the volume element corresponds precisely 
to finding a skew-symmetric phase-space metric (symplectic form)
that is an integral of motion. 
Therefore we consider the skew-symmetric metric. 

Sergi \cite{Sergi2,Sergi3} uses the skew-symmetric phase-space metric  
that is not satisfied by the Jacoby identity.  
As the result the generalization of the Poisson 
brackets for non-Hamiltonian systems leads one to non-Lie algebra.
Note that non-Lie algebra for non-Hamiltonian systems is considered 
in \cite{Tartmf4}. 
In this paper we consider an invariant antisymmetric 
phase-space metric that satisfies the Jacoby identity, and 
defines the Lie algebra in phase-space.
We call the metric is invariant if the metric tensor field 
is an integral of motion ($d\omega/dt=0$). 
This invariant phase-space metric $\omega_{kl}({\bf x},t)$ defines
the invariant phase-space measure $v$ by the equation 
$v=(1/n!)\omega^n=\sqrt{g({\bf x},t)} dx^1\wedge...\wedge dx^{2n}$,
where $g({\bf x},t)$ is the metric determinant
$g({\bf x},t)=det(\omega_{kl}({\bf x},t))$ and $dv/dt=0$. 
The suggested time-dependent skew-symmetric phase-space 
metric leads to a constant value of the entropy density, so that the 
associated distribution function obeys an evolution equation 
associated with incompressible dynamical flow. 

Note that the invariant phase-space metric 
of some non-Hamiltonian systems can leads us to the lack of smoothness 
of the metric. In this case, the phase-space probability distribution 
can be collapsed onto a fractal set of dimensionality lower 
than in the Hamiltonian case \cite{Dorfman,HPAK}.
Unfortunately the description of lack of smoothness 
in \cite{Dorfman,HPAK}
is considered without using the curved phase-space 
approach \cite{Tuck1,Tuck2,Tuck3}.
Note that  classical systems that are Hamiltonian systems 
in the usual phase-space
are non-Hamiltonian systems in the fractional phase-space \cite{chaos,PRE}.

In the papers \cite{Tarpla1,Tarmsu}, the quantization
of the evolution equations for non-Hamiltonian and dissipative 
systems was suggested.
Using this quantization it is easy to derive the quantum
analog of the invariant Poisson brackets, which satisfy 
the rule of term-by-term differentiation with respect to time.


\end{document}